\documentclass{amsart}

\usepackage{amssymb}
\usepackage{amsmath}
\usepackage{bm}

\newcommand{\half}{\frac{1}{2}}
\newcommand{\thalf}{\tfrac{1}{2}}

\newcommand{\sums}{\mathop{{\sum}^{*}}}

\newcommand{\sym}{{\rm sym}}
\newcommand{\intt}{\int_{-\infty}^{\infty}}

\numberwithin{equation}{section}

\newtheorem{theorem}{Theorem}[section]

\newtheorem{lemma}[theorem]{Lemma}

\makeatletter
\@namedef{subjclassname@2010}{%
  \textup{2010} Mathematics Subject Classification}
\makeatother

\begin{document}

\title{On subconvexity bounds for twisted $L$-functions}

\author{Rizwanur Khan}
\address{
Department of Mathematics\\ University of Mississippi\\ University, MS 38677}
\email{rrkhan@olemiss.edu }

\thanks{The author was supported by the National Science Foundation Grant DMS-2001183 and the Simons Foundation (award 630985). Any opinions, findings and conclusions or recommendations expressed in this material are those of the authors and do not necessarily reflect the views of the National Science Foundation.}

\begin{abstract} 
We prove hybrid subconvexity bounds for a wide class of twisted $L$-functions $L(s,f\times \chi)$ at the central point, including a new instance of the Weyl subconvexity bound.
\end{abstract}

\maketitle

\section{Introduction}

The subconvexity problem is an important problem in the theory of $L$-functions concerned with bounding an automorphic $L$-function $L(s)$ on the critical line $\Re(s)=\half$. A bound which can readily be obtained using essentially only the functional equation of the $L$-function is called the convexity bound and is given by
\begin{align}
\label{conv} L(\thalf + it) \ll \mathfrak{q}^{\kappa+\epsilon},
\end{align}
where $\mathfrak{q}$ is the analytic conductor of the $L$-function \cite[Section 5.1]{iwakow}, $\kappa=\frac14$ and $\epsilon>0$ is arbitrarily small. The subconvexity problem is to establish (\ref{conv}) with some number $\kappa$ strictly less than $\frac14$. The best possible expectation, essentially, is that we may take $\kappa=0$, and this is the famous Lindel\"{o}f conjecture. The subconvexity bound for $GL(1)$ $L$-functions is classical, and for $GL(2)$ automorphic $L$-functions it was completely resolved by Michel and Venkatesh \cite{micven}. However it is of great interest to not only break the convexity bound, but to reach certain milestones called the Burgess exponent ($\kappa=\frac{3}{16}$) and the Weyl exponent ($\kappa=\frac{1}{6}$). These are named after classical results in $GL(1)$, but also occur in higher rank. It is unknown why the Burgess and Weyl exponents occur in so many unrelated situations, but in all cases they seem to represent natural barriers which can require deep ideas to overcome. 

 We mention a few examples in $GL(1)$ and $GL(2)$ at the central point $s=\half$. For Dirichlet $L$-functions of modulus $q$, the analytic conductor is equal to $q$. Burgess \cite{bur} proved the bound $
L(\half ,\chi) \ll q^{\frac{3}{16}+\epsilon}$, Conrey and Iwaniec \cite{coniwa} proved the Weyl bound $q^{\frac{1}{6}+\epsilon}$ for real characters,  Petrow and Young \cite{petyou1, petyou2} proved the Weyl bound for all Dirichlet characters, and Mili\'{c}evi\'{c} \cite{mil} proved a sub-Weyl bound in the case of Dirichlet characters of prime power modulus $p^n$ with $n$ large. Turning now to $GL(2)$, let $f$ be a holomorphic Hecke cusp form or a Hecke Maass cusp form for the full modular group, with associated $L$-function $L(f,s)$. The Weyl subconvexity bound was proven for $L(\half,f)$ in the spectral and weight aspects by Ivi\'{c} \cite{ivi} and Peng \cite{pen} respectively. 
If $f$ is a newform of level $q$ and trivial nebentypus, the analytic conductor of $L(s,f)$ in the level aspect is $q$. Duke, Friedlander and Iwaniec \cite{dfi} proved the first subconvexity bound $L(\half, f)\ll q^{1-\frac{1}{192}+\epsilon}$. This has seen some recent improvements. The current best result for $q$ an odd prime is $\kappa=\frac{1}{5}$ conditionally on the Ramanujan conjecture, by Kiral and Young \cite{kiryou} and Blomer and Khan \cite{blokha2}. Unconditionally it is $\kappa=0.2073$ by Blomer, Humphries, Khan, Milinovich \cite{bhkm}. For newforms with primitive nebentypus, the best result is due to Blomer and Khan \cite{blokha}, with exponent $\kappa=\frac14-\frac{1}{128}$.

It is highly frustrating that the subconvexity exponent in the level aspect does not attain even the Burgess milestone, despite better results in other $GL(2)$ aspects. There are however exceptions for special forms which arise as twists. 
%For example, if $f$ is of level $1$ and we twist by a primitive Dirichlet character $\chi$ of level $q$, then $f\times \chi$ is a newform of level $q^2$ and nebentypus $\chi^2$.
For twists of a form $f$ of {\it fixed} level by a primitive Dirichlet character $\chi$ of modulus $q$, the Burgess bound for $L(\half, f\times \chi)$ in the $q$ aspect is known by the results of Bykovski\u\i  \ \cite{byk} and Blomer, Harcos, and Michel \cite{bloharmic, blohar}, and the Weyl bound is known when $\chi$ is real by work of Petrow and Young \cite{petyou3}. The Weyl bound was also established for twists by characters of prime power moduli by Blomer and Mili\'{c}evi\'{c} \cite{blomil} and Munshi and Singh \cite{munsin}. Outstanding results also exist for twists of forms whose level is not fixed. Conrey and Iwaniec \cite{coniwa} proved the Weyl bound for $L(1/2,f\times \chi)$ for newforms $f$ of level $q$ and trivial nebentypus, twisted by real primitive characters $\chi$ of modulus $q$ (see also \cite{you}). Petrow and Young \cite{petyou1} generalized this result to allow any primitive character $\chi$ of modulus $q$ as long as $f$ is a newform of level $q$ and nebentypus $\overline{\chi}^2$. The works of Iwaniec and Conrey and of Petrow and Young are particularly significant because they also lead to the Weyl subconvexity bound for Dirichlet $L$-functions.

The goal of this paper is to prove subconvexity bounds for a wide class of twisted $L$-functions, including a new instance of the Weyl subconvexity bound. We restrict to prime level and modulus in order to minimize technical details and concentrate more on the main ideas.
\begin{theorem}\label{main}
Let $\epsilon>0$. Let $f$ be a holomorphic newform with weight $k\ge 2$, prime level $q$ and trivial nebentypus. Let $\chi$ be a primitive Dirichlet character of prime modulus $p$. Suppose that $p^\epsilon\le q\le p^{2+\epsilon}$ and $(q,p)=1$. We have 
\begin{align*}
L(\thalf, f\times \chi)\ll_{k,\epsilon} p^\epsilon (q^{\half}+p^{\half}).
\end{align*}
\end{theorem}
\noindent By Atkin-Lehner theory \cite[Proposition 14.20]{iwakow}, $L(s, f\times \chi)$ is the $L$-function attached to a newform of level $qp^2$ and nebentypus $\chi^2$. Thus the result above gives the convexity bound or better in all cases, a subconvexity bound when $q<p^{2-\delta}$ for some $\delta>0$, the Burgess subconvexity bound or better for $p^{\frac23-\epsilon}<q<p^{\frac65+\epsilon}$, and the Weyl subconvexity bound for $p^{1-\epsilon}<q<p^{1+\epsilon}$, which is the best our result can do. This last bound may be viewed as complementary to the work of the Conrey and Iwaniec \cite{coniwa} and Petrow and Young \cite{petyou1}. They obtained the Weyl bound when the level of the form and modulus of the Dirichlet character coincide, while we do the same when the level and modulus are of the same size, but coprime.

Theorem \ref{main} is a type of hybrid subconvexity bound, since it yields subconvexity in both parameters $p$ and $q$ simultaneously. Such hybrid bounds are of great interest and have been studied by Blomer and Harcos \cite{blohar, blohar2}, Aggarwal, Jo, and Nowland \cite{aggjonow}, and Chen and Hou \cite{houche}. Of these results, the last one yields subconvexity for the widest range: for $q$ as large as $p^{\frac32-\delta}$, but it falls short of any milestone subconvexity bound. Our result works in an even wider range, and gives much stronger bounds when $q$ and $p$ are about the same size. Moreover our proof is completely different and relatively simple. In particular, it totally avoids the shifted convolution problem, which was central in the aforementioned works on hybrid subconvexity. Theorem \ref{main} follows immediately by establishing the following second moment estimate.

\begin{theorem} \label{2ndmoment}  Keep the notation in the statement of Theorem \ref{main}. Let $B_k^\star(q)$ denote a basis of holomorphic newforms of level $q$ and weight $k\ge 2$. We have
\begin{align} 
\label{2ndm}  \sum_{f\in B_k^*(q)} |L(\thalf, f\times\chi)|^2  \ll_{k,\epsilon}  p^{\epsilon}(q+p)
\end{align}
for any $\epsilon>0$. 
\end{theorem}
\noindent For comparision, the large sieve \cite[Theorem 7.26]{iwakow} would give a bound of $O(p^\epsilon(q+q^\half p))$ for the second moment. In connection to this, we mention a further application that our work may have. Liu, Masri, and Young \cite[page 16]{liumasyou} observed that if \cite[Theorem 1.2]{butkha} could be extended to small weights, then the large sieve estimate for the left hand side of \eqref{2ndm} would imply a hybrid subconvexity bound for $L(\half, f\times \Theta_\chi)$ that is uniform in {\it all} ranges of $q$ and $p$. Here $\Theta_\chi$ is the theta series associated to an ideal class group character of $\mathbb{Q}(\sqrt{-p})$ . We note that using our Theorem \ref{2ndmoment} instead of the large sieve would give a stronger subconvexity bound for this $GL(2)\times GL(2)$ $L$-function.

Throughout, we follow the $\epsilon$-convention: that is, $\epsilon$ will always be positive number which can be taken as small as we like, but may differ from one occurence to another. All implied constants may depend on $\epsilon$ and $k$.

{\bf Acknowledgement.} I am grateful to Matthew P. Young for his interest in this work and many helpful comments.

\section{Preliminaries}

\subsection{Poisson summation}

One of the basic tools we will use is the Poisson summation formula. We state it here.
\begin{lemma}
Let $\psi$ be a function in the Schwartz class. We have
\[
\sum_{\substack{n\in \mathbb{Z}\\ n\equiv a \bmod r}} \psi(n) = \frac{1}{r} \sum_{\substack{m\in \mathbb{Z}}} e\Big(\frac{am}{r}\Big) \hat{\psi}\Big(\frac{m}{q} \Big),
\] 
where $\hat{\psi}(\xi)=\intt \psi(t) e(-t\xi) dt$  and $e(x)=e^{2\pi i x}$.
\end{lemma}

\subsection{Modular forms and $L$-functions}
Let $S_k(q)$ denote the space of holomorphic cusp forms of prime level $q$, weight $k$, and trivial nebentypus. Let $S_k^*(q)\subset S_k(q) $ denote the space of newforms. Every $f\in S_k(q)$ has a Fourier series expansion
\begin{align*}
f(z)=\sum_{n=1}^\infty a_f(n) n^{\frac{k-1}{2}} e(nz)
\end{align*}
for $\Im(z)>0$. Let $B_k(q)$ denote an orthogonal basis of $S_k(q)$ which contains a basis $B_k^\star(q)$ of $S_k^*(q)$, normalized so that $a_f(1)=1$ for every $f\in B_k(q)$.

Let $\chi$ be primitive Dirichlet character of prime modulus $p$. For $f\in B_k(q)$, define
\begin{align*}
(f\times \chi)(z)=\sum_{n=1}^\infty \chi(n) a_f(n) n^{\frac{k-1}{2}} e(nz).
\end{align*}
For $f\in B_k^\star(q)$, this is a newform of level $qp^2$, weight $k$, and nebentypus $\chi^2$. The associated $L$-function is entire and for $\Re(s)>1$ equals
\begin{align*}
L(s,f\times \chi)=\sum_{n=1}^\infty \frac{\chi(n) a_f(n)}{ n^s}.
\end{align*}
This satisfies the functional equation
\begin{align}
\label{feq} \Lambda(s, f\times \chi):= \Big(\frac{p q^\half  }{2\pi}\Big)^s\Gamma(s+\tfrac{k-1}{2}) L(s,f\times \chi)= \varepsilon \Lambda(1-s, f\times \overline{\chi}),
\end{align}
for some complex number $\varepsilon$ of modulus 1, which depends on $q, p,$ and $\chi$. Thus the analytic conductor at $s=\frac12$ is $kp^2q$. These facts may be found in \cite[Section 14]{iwakow}. Let $L(s,\sym^2 f)$ denote the symmetric square  attached to $f\in B_k^*(q)$. This is defined in \cite[Section 5.12]{iwakow}, but all we need are the bounds
\begin{align}
\label{sym1} (kq)^{-\epsilon}\ll  L(1,\sym^2 f) \ll  (kq)^{\epsilon}
\end{align}
given in \cite{hofloc}. In particular, $L(1,\sym^2 f)$ is positive.

\subsection{Approximate functional equation}

We will need the following consequence of the approximate functional equation \cite[Theorem 5.3]{iwakow}.
\begin{lemma} \label{afe}
\begin{align}
\label{afe-line} |L(\thalf,f\times \chi)|^2\ll \Big| \sum_{n\ge 1 } \frac{a_f(n) \chi(n) }{ n^\half}V\Big(\frac{n}{q^\half p}\Big)\Big|^2+\Big| \sum_{n\ge 1 } \frac{a_f(n) \overline{\chi}(n) }{ n^\half}V\Big(\frac{n}{q^\half p}\Big)\Big|^2,
\end{align}
where we define the real function
\begin{align*}
V(x)=\frac{1}{2\pi i} \int_{(2)} (2\pi x)^{-u}  \frac{\Gamma(\frac{k}{2}+u)}{\Gamma(\frac{k}{2})} \frac{du}{u}
\end{align*}
for $x>0$.
\end{lemma}
\noindent By a standard argument of shifting contours, we may restrict the sum above to $n\le q^{\half} p^{1+\epsilon}$, up to an error of $O(p^{-100}$). Further, by splitting the $n$-sum in (\ref{afe-line}) into $O(\log q)$ sums over dyadic intervals and applying the Cauchy-Schwarz inequality, for the purposes of Theorem \ref{2ndmoment} it suffices to prove
\begin{align}
 \label{toprove} \frac{1}{N} \sum_{f\in B_k^*(q)} \Big| \sum_{n\ge 1} \chi(n)  a_f(n) V\Big(\frac{n}{N}\Big)\Big|^2  \ll  p^\epsilon(q+p)
\end{align}
for any
\[
1\le N \le q^\half p^{1+ \epsilon}
\] 
and any smooth function $V$ compactly supported on the positive reals with $\|V^{(j)}\|_\infty \ll_j (p^\epsilon)^{j+1}$ for all $j\ge 0$. 

\subsection{Petersson trace formula} 

The Petersson trace formula \cite[Proposition 2.1]{iwamic} states that
\begin{align}
\label{trace} \sum_{f\in B_k(q)} \frac{\Gamma(k-1)}{(4\pi)^{k-1}\langle f, f\rangle } a_f(n_1)\overline{a_f(n_2)}=\delta_{(n_1=n_2)}+2\pi  i^{-k}\sum_{c\ge 1} \frac{S(n_1,n_2,cq)}{cq} J_{k-1}\Big(\frac{4\pi \sqrt{n_1n_2}}{cq}\Big),
\end{align}
where $\langle f , f \rangle$ is the Petersson inner product, $J_{k-1}(x)$ is a $J$-Bessel function, and $\delta_P$ equals $1$ if the statement $P$ is true, and $0$ otherwise.

By \cite[equations (2.3, 2.24)]{iwamic}, we have
\begin{align*}
\frac{2\pi^2}{q(k-1)L(1,\sym^2 f)}=\frac{\Gamma(k-1)}{(4\pi)^{k-1}\langle f, f\rangle }
\end{align*}
for $f\in  B^*_k(q)$. This combined with the bound \eqref{sym1} means that for \eqref{toprove}, it suffices to prove 
\begin{align*}
 \frac{1}{N} \sum_{f\in B^*_k(q)}  \frac{\Gamma(k-1)}{(4\pi)^{k-1}\langle f, f\rangle }  \Big| \sum_{n\ge 1}  \chi(n)  a_f(n)V\Big(\frac{n}{N}\Big)\Big|^2  \ll  p^\epsilon(q+p).
\end{align*}
By positivity, we may enlarge $B^*_k(q)$ to $B_k(q)$, expand out the square $|\displaystyle \sum_{n\ge 1} \ldots|^2$, and apply the Petersson trace formula \eqref{trace}. After doing so, we need to prove
\begin{align}
\label{toprove2} \Big|\frac{1}{N} \sum_{\substack{n_1,n_2\ge 1\\ n_1=n_2}} V\Big(\frac{n_1}{N}\Big) V\Big(\frac{n_2}{N}\Big)\Big| + \Big|\frac{1}{N}\sum_{n_1,n_2,c\ge 1} \frac{S(n_1,n_2,cq) \chi(n_1)  \overline{\chi}(n_2)  }{cq}  V\Big(\frac{n_1}{N}\Big) V\Big(\frac{n_2}{N}\Big) J_{k-1}\Big(\frac{\sqrt{n_1n_2}}{cq}\Big)\Big|   \ll  p^\epsilon\Big(1+\frac{p}{q}\Big).
\end{align}
The contribution of the `diagonal' $\sum_{n_1=n_2}$ is obviously $O(p^\epsilon)$. For the off-diagonal we first note the following.

\begin{lemma}\label{Bx}
For $x\le p^{\epsilon}$ we have 
\[
J_{k-1}(x)=xW_1(x)
\] where $W_1(x)$ is a smooth function satisfying $\|W_1^{(j)}\|_\infty \ll_j (p^\epsilon)^{j+1}$ for all $j\ge 0$ and $x$ in this range.

For $x> p^{\epsilon}$, we have
\begin{align}
\label{bx-main} J_{k-1}(x)=\Im\Big( \frac{e(2x)}{\sqrt{x}}W_2(x) \Big),
\end{align}
where $\Im$ denotes the imaginary part and $W_2$ is a smooth function satisfying
\begin{align*}
 x^jW_2^{(j)}(x)\ll_{j,\epsilon} (p^{\epsilon})^{j+1}
\end{align*}
for all $j\ge 0$ and $x$ in this range. 
\end{lemma}

\proof 

When $x\le p^{-\epsilon}$, use the power series \cite[page 82]{iwa},
\[
J_{k-1}(x)=\sum_{\ell=0}^\infty \frac{(-1)^\ell (x/2)^{k-1+2\ell}}{\ell! (k+\ell)!}.
\]
When $x>p^{-\epsilon}$, use \cite[section 4]{iwamic}.
\endproof 

\subsection{Conclusion of the set-up}
We now apply Lemma \ref{Bx} to \eqref{toprove2}. The extra factors $\frac{\sqrt{n_1 n_2}}{cq}$ or $(\frac{\sqrt{n_1 n_2}}{cq})^{-1/2}$ can be written as $\frac{N}{cq}$ and $(\frac{N}{cq})^{-1/2}$ after redefining $V(\frac{n_1}{N})$ and $V(\frac{n_2}{N})$. By considering the real and imaginary parts of $W_1$ and $W_2$ separately, we can assume that they are real functions. Thus it suffices to prove, for each sign $\pm$, that
\[
\mathcal{S}_1+\mathcal{S}_2^\pm\ll p^\epsilon\Big(1+\frac{p}{q}\Big),
\]
where
\begin{align}
\label{s1def} &\mathcal{S}_1=\frac{1}{ C_1^2 q^2}\sum_{C_1\le c\le 2C_1}\Big| \sum_{\substack{n_1,n_2\ge 1 }} S(n_1,n_2,cq) \chi(n_1)  \overline{\chi}(n_2)  V\Big(\frac{n_1}{N},\frac{n_2}{N}\Big)\Big|,\\
\label{s2def} &\mathcal{S}_2^\pm=\frac{1}{N^\frac{3}{2} q^\half C_2^\half}\sum_{C_2\le c\le 2C_2}\Big| \sum_{\substack{n_1,n_2\ge 1 }}  S(n_1,n_2,cq) \chi(n_1)  \overline{\chi}(n_2)  e\Big(\frac{\pm 2\sqrt{n_1n_2}}{cq}\Big)  V\Big(\frac{n_1}{N},\frac{n_2}{N}\Big)\Big|,
\end{align}
for any 
\begin{align}
\label{crange} C_1\ge \frac{N}{ p^{\epsilon} q} , \ \ \ \ \ \  1\le C_2 \le \frac{N}{ p^{\epsilon} q},
\end{align}
and where we set
\[
V(x_1,x_2):= V_i(x_1,x_2)=V(x_1)V(x_2)W_i\Big(\frac{N\sqrt{x_1x_2}}{cq}\Big) 
\]
for $i=1,2$. We drop the subscript from $V$ because all we need is that $V$ is compactly supported on $\mathbb{R}^{+}\times \mathbb{R}^{+}$ and satisfies
 \[
\frac{\partial^{i}}{\partial x_1^i } \frac{\partial^{j}}{\partial x_2^j } V(x_1,x_2)  \ll_j (p^\epsilon)^{1+i+j}
\] for all $i,j\ge 0$.

\section{Proof of Theorem \ref{2ndmoment}}

The goal now is to prove the required bounds for $\mathcal{S}_1$ and $\mathcal{S}_2^\pm$. 

\subsection{Possion summation in $n_1$} By separating $n_1$ into residue classes modulo $cpq$, and applying Poisson summation, we get that

\begin{align}
\label{poi1} &\mathcal{S}_1=   \frac{N}{C_1^3 q^3 p}\sum_{C_1\le c\le 2C_1}\Big|  \sum_{\substack{n_2\ge 1 }} \sum_{m_1\in \mathbb{Z}}  \  \sum_{a
\bmod cpq} S(a,n_2,cq) \chi(a)  \overline{\chi}(n_2)  e\Big(\frac{am_1}{cpq}\Big) \mathcal{V}_1\Big( m_1, \frac{n_2}{N}\Big)\Big|, \\
\label{poi2} & \mathcal{S}_2^\pm=\frac{1}{N^\half C_2^\frac32 q^\frac32 p}\sum_{C_2\le c\le 2C_2 }\Big|\sum_{\substack{n_2\ge 1}} \sum_{m_1\in \mathbb{Z}}  \  \sum_{a
\bmod cpq} S(a,n_2,cq) \chi(a)  \overline{\chi}(n_2)  e\Big(\frac{am_1}{cpq}\Big) \mathcal{V}_2^\pm\Big( m_1, \frac{n_2}{N}\Big)\Big|,
\end{align}
where
\begin{align}
\label{vint1} &\mathcal{V}_1(m_1,y)= \intt e\Big(-\frac{xm_1 N}{cpq}\Big)V(x,y)dx.\\
\label{vint2} &\mathcal{V}_2^\pm(m_1,y)= \intt e\Big(\frac{\pm 2\sqrt{x y N^2}}{cq}\Big)  e\Big(-\frac{xm_1 N}{cpq}\Big)V(x,y)dx.
\end{align}

\begin{lemma}\label{statp}
Suppose $c\asymp C_1$ in $\mathcal{V}_1(m_1,y)$ and $c\asymp C_2$ in $\mathcal{V}_2(m_1,y)$. Suppose $m_1\neq 0$. We have
\begin{align*}
&\mathcal{V}_1(m_1,y)= W(y) \delta_{(|m_1|\le  \frac{C_1qp^{1+\epsilon}}{N})}+O(p^{-100}),\\
&\mathcal{V}_2^{+}(m_1,y)= \Big(\frac{cq}{N}\Big)^\half e\Big(\frac{ypN}{cqm_1}\Big)W(y)\delta_{(p^{1-\epsilon}\le m_1\le p^{1+\epsilon})}+O(p^{-100}),\\
&\mathcal{V}_2^{-}(m_1,y)= \Big(\frac{cq}{N}\Big)^\half e\Big(\frac{ypN}{cqm_1}\Big)W(y)\delta_{(-p^{1+\epsilon}\le m_1\le -p^{1-\epsilon})}+O(p^{-100}),\\
\end{align*}
where $W(y)$ is a smooth compactly supported function on the positive reals, depending on $\pm, c,p,q, N, m_1$, with $\|W^{(j)}\|_\infty \ll_j (p^\epsilon)^j$ for all $j\ge 0$.

\end{lemma}
\proof

For $\mathcal{V}_1(m_1,y)$, we are simplying renaming this function to $W(y)$. It's clear from the expression \eqref{vint1} that $W$ satisfies the required properties. Further, by integrating by parts multiple times in \eqref{vint1}, integrating $e(-\frac{xm_1N}{cpq})$ and differentiating the rest of the integrand, we get that \eqref{vint1} is $O(p^{-100})$ unless 
\[
|m_1|\le \frac{C_1 q p^{1+\epsilon}}{N}.
\]

For $\mathcal{V}_2^\pm (m_1,y)$, we apply a stationary phase approximation. Let
\[
h\pm(x):= 2\pi\Big( \frac{\pm 2\sqrt{x y N^2}}{cq}-\frac{xm_1 N}{cpq}\Big)
\]
be the phase of the integrand. 

When $|m_1|< p^{1-\epsilon}$ or $|m_1|>p^{1+\epsilon}$, we have $|h'(x)|\gg \frac{N}{cq} \gg p^\epsilon$. We can apply \cite[Lemma 8.1]{bky} with the parameters $R= Y=\frac{N}{cq}$, $X=Q=U=p^\epsilon$ to see that the integral in \eqref{vint2} is $O(p^{-100})$. 

Now suppose $p^{1-\epsilon}\le |m_1| \le p^{1+\epsilon}.$ The stationary point $x_0$ satisfies $h'(x_0)=0$, which implies
\begin{align}
\label{statpt} m_1 \sqrt{x_0}= \pm p \sqrt{y}.
\end{align}
 In the $+$ case, there is a unique solution when $p^{1-\epsilon}\le  m_1\le p^{1+\epsilon}$.  In $-$ case, there is a unique solution when $-p^{1+\epsilon}\le m_1\le -p^{1-\epsilon}$. We show the details for $\mathcal{V}_2^+$ and $m_1>0$ only. From \eqref{statpt}, we get that
\[
x_0=\frac{p^2 y}{m_1^2}.
\]
We have
\[
h(x_0)= 2\pi  \frac{ypN}{cqm_1}
\]
and
\[
h''(x_0)=\frac{-m_1^3}{2yp^3} \frac{N}{cq}.
\]
Applying \cite[Proposition 8.2]{bky} with the parameters $Q=1$, $X=V=V_1=p^\epsilon$, $Y=\frac{N}{cq}$, $\delta=\frac{1}{9}$, we get that
\[
\mathcal{V}_{+}(m_1,y)=\Big(\frac{cq}{N}\Big)^\half e\Big(\frac{ypN}{cqm_1}\Big)W(y)\delta_{(p^{1-\epsilon}\le m_1\le p^{1+\epsilon})}+O(p^{-100}),
\]
where $W(y)$ is as described in the statement of the lemma. The $\delta$ function indicates the range of $m_1$, but anyway this information is contained in $W(y)$.
\endproof

We now evaluate the arithmetic part in \eqref{poi1} and \eqref{poi2}. We first observe that we can restrict to $(c,p)=1$ in \eqref{s1def} and \eqref{s2def}, up to an admissible error. This is because in \eqref{s2def} we have $c\ll \frac{N}{p^\epsilon q} \ll \frac{p^{1+\epsilon}}{q^\half}$ and so $c$ is too small to be divisible by $p$. In \eqref{s1def}, we  can trivially bound the terms with $p|c$, using Weil's bound for the Kloosterman sum, to see that their contribution is  $O(\frac{p^{1/2+\epsilon}}{q^{1/2}})$, which is dominated by the required bound in \eqref{toprove2}. Having restricted to $(c,p)=1$, we can write the complete sum mod $cpq$ in terms of complete sums mod $cq$ and mod $p$. We get, writing $\sums$ to denote a sum over primitive residue classes and a bar over a residue class to mean the inverse of that residue class, 
\begin{align}
\nonumber \sum_{a
\bmod cpq} S(a,n_2,cq) \chi(a)   e\Big(\frac{am_1}{cpq}\Big) &=  \sum_{\substack{a_1 \bmod cq\\ a_2 
\bmod p}} S(a_1,n_2,cq) \chi(a_2)   e\Big(\frac{a_1m_1 \overline{p} }{cq}\Big)e\Big(\frac{a_2 m_1 \overline{cq} }{p}\Big)\\
\nonumber &=  \sum_{\substack{a_1 \bmod cq\\ a_2 
\bmod p}} \ \sums_{b\bmod cq}  e\Big(\frac{a_1\overline{b} + n_2 b}{cq}\Big) \chi(a_2)   e\Big(\frac{a_1m_1 \overline{p} }{cq}\Big)e\Big(\frac{a_2 m_1 \overline{cq} }{p}\Big)\\
\label{aritheval} &=  cq \tau(\chi) \delta_{(m_1,cq)=1}   e\Big(\frac{-n_2 p \overline{m_1} }{cq}\Big) \chi(m_1) \overline{\chi}(cq),
\end{align}
where $\tau(\chi)$ is the Gauss sum. The final equality uses \cite[equation (3.12)]{iwakow} and the orthogonality of additive characters \cite[page 44]{iwakow}.

\

After applying \eqref{aritheval} and Lemma \ref{statp} to \eqref{poi1} and \eqref{poi2},  and using that $|\tau(\chi)|=p^\half$, we see that it suffices to prove 
\begin{align}
\label{p1a} & \frac{N}{C_1^2 q^2  p^\frac12} \sum_{\substack{C_1 \le c\le 2C_1 }} \ \sum_{\substack{|m_1|\le \frac{C_1qp^{1+\epsilon}}{N}\\ (m_1,cq)=1 }}  \Bigg|   \sum_{n_2\ge 1}  e\Big(\frac{-n_2 p \overline{m_1} }{cq}\Big) \overline{\chi}(n_2) W\Big(\frac{n_2}{N}\Big) \Bigg| \ll p^\epsilon\Big(1+\frac{p}{q}\Big), \\
\label{p2a} & \frac{1}{N  p^\frac12} \sum_{\substack{C_1 \le c\le 2C_2 }} \ \sum_{\substack{p^{1-\epsilon}\le |m_1| \le p^{1+\epsilon} \\ (m_1,cq)=1  }} \Bigg|   \sum_{n_2\ge 1}  e\Big(\frac{-n_2 p \overline{m_1} }{cq}\Big) e\Big(\frac{n_2 p}{cqm_1}\Big) \overline{\chi}(n_2) W\Big(\frac{n_2}{N}\Big) \Bigg|\ll p^\epsilon\Big(1+\frac{p}{q}\Big).
\end{align}

\subsection{Reciprocity}

Using reciprocity, or the Chinese remainder theorem, for $(m_1,cq)=1$, we have
\[
e\Big(\frac{n_2 p}{cqm_1}\Big) = e\Big(\frac{n_2 p\overline{m_1}}{cq}\Big) e\Big(\frac{n_2 p\overline{cq} }{m_1}\Big).
\]
Thus \eqref{p1a} and \eqref{p2a} become 
\begin{align}
\label{p1b} & \ \frac{N}{C_1^2 q^2  p^\frac12} \sum_{\substack{C_1 \le c\le 2C_1 }} \ \sum_{\substack{|m_1|\le \frac{C_1qp^{1+\epsilon}}{N}\\ (m_1,cq)=1 }}   \Bigg|   \sum_{n_2\ge 1}  e\Big(\frac{n_2 p\overline{cq} }{m_1}\Big) e\Big(\frac{-n_2 p}{cqm_1}\Big)  \overline{\chi}(n_2) W\Big(\frac{n_2}{N}\Big) \Bigg| \ll p^\epsilon\Big(1+\frac{p}{q}\Big), \\
\label{p2b} & \frac{1}{N  p^\frac12} \sum_{\substack{C_1 \le c\le 2C_2 }} \ \sum_{\substack{p^{1-\epsilon}\le |m_1| \le p^{1+\epsilon} \\ (m_1,cq)=1  }} \Bigg|   \sum_{n_2\ge 1}  e\Big(\frac{n_2 p\overline{cq} }{m_1}\Big) \overline{\chi}(n_2) W\Big(\frac{n_2}{N}\Big)  \Bigg|\ll p^\epsilon\Big(1+\frac{p}{q}\Big).
\end{align}

\subsection{Poisson summation in $n_2$}

\subsubsection{The case of large $c$}\label{ss} We first consider the left hand side \eqref{p1b}. Observe that we can add the condition $(m_1,p
)=1$ to the sum, because the total contribution of the terms with $p|m_1$ is trivially $O(\frac{p^{1/2+\epsilon}}{q^{1/2}})$. The following notation will assume that $m_1$ is positive, with the case of negative $m_1$ being entirely similar. Define 
\begin{align}
\label{wint} \mathcal{W}(m_2)= \intt   e\Big(\frac{-ym_2 N}{m_1 p}\Big) e\Big(\frac{-yN p}{cqm_1}\Big) W(y) dy,
\end{align}
and keep in mind that $W$ is compactly supported.
By separating $n_2$ into residue classes modulo $m_1p$ and appying Poisson summation, we get that 
\begin{align}
\nonumber &\sum_{n_2\ge 1}  e\Big(\frac{n_2 p\overline{cq} }{m_1}\Big)  \overline{\chi}(n_2)  e\Big(\frac{-n_2 p}{cqm_1}\Big) W\Big(\frac{n_2}{N}\Big)\\
\nonumber &= \frac{N}{m_1 p}\sum_{b\bmod m_1p}  \ \sum_{m_2\in \mathbb{Z}} e\Big(\frac{b p\overline{cq} }{m_1}\Big)  \overline{\chi}(b) e\Big(\frac{bm_2}{m_1p}\Big) \mathcal{W}(m_2)\\
\nonumber &= \frac{N}{m_1 p}\sum_{\substack{b_1 \bmod m_1\\ b_2 \bmod p}}  \ \sum_{m_2\in \mathbb{Z}} e\Big(\frac{b_1 p\overline{cq} }{m_1}\Big)  \overline{\chi}(b_2) e\Big(\frac{b_1m_2\overline{p}}{m_1}\Big) e\Big(\frac{b_2m_2\overline{m_1}}{p}\Big)  \mathcal{W}(m_2)\\
\label{arp2} &= \frac{N \tau(\overline{\chi})  }{p} \sum_{\substack{m_2\in \mathbb{Z}\\ cm_2\equiv p^2\overline{q}\bmod m_1 }}  \overline{\chi}(m_2)\chi(m_1) \mathcal{W}(m_2).
\end{align}
By integrating by parts multiple times in \eqref{wint}, we find that $\mathcal{W}(m_2)\ll p^{-100}$ unless
\[
\Big|\frac{m_2N}{m_1p}+\frac{Np}{cqm_1}\Big|\le p^\epsilon,
\]
in which case $\mathcal{W}(m_2)\ll 1$. This implies
\[
c m_2=\Big\lfloor \frac{-p^2}{q} \Big\rfloor+\ell
\]
for
\[
|\ell|\le \frac{c m_1 p^{1+\epsilon} }{N}.
\]
Since $\ell$ has only $O(\ell^\epsilon)$ divisors, once we apply \eqref{arp2} to \eqref{p1b} and bound everything absolutely, we may collapse the the sum over $c$ and $m_2$ to a sum over $\ell$ only. The congruence condition in \eqref{arp2} translates to a congruence condition on $\ell$. Thus we see that to establish \eqref{p1b}, it suffices to prove that for any $\alpha_{m_1}$, we have
\[
\frac{N^2}{C_1^2 q^2 p}  \sum_{1\le m_1\le \frac{C_1pq^{1+\epsilon}}{N}}\ \sum_{\substack{  |\ell| \le  \frac{C_1 m_1 p^{1+\epsilon} }{N} \\ \ell \equiv \alpha_{m_1} \bmod m_1 }} 1 \ll p^\epsilon\Big(\frac{p}{q}+1\Big).
\]
This follows since the $\ell$-sum is $O(\frac{C_1 p^{1+\epsilon} }{N}+1).$

\subsubsection{The case of small $c$} We now consider the left hand side \eqref{p2b}, assuming $m_1>0$. Again, a trivial bound for the terms with $p|m_1$ allows us to restrict to $(m_1,p)=1$. The details of Poisson summation in $n_2$ are the same as in subsection \ref{ss}, except that the transform function will equal
\[
\mathcal{W}(m_2)= \intt   e\Big(\frac{-ym_2 N}{m_1 p}\Big)  W(y) dy,
\]
from which it follows using integration by parts that we can restrict to
\[
|m_2|\le p^\epsilon \frac{m_1p}{N}\le \frac{  p^{2+\epsilon}}{N}.
\]
The same congruence condition as in \eqref{arp2} will hold here. Thus to establish \eqref{p2b}, it suffices to prove that for any $\alpha_{m_1}$,
\[
\frac{1}{p}  \sum_{\substack{p^{1-\epsilon}\le m_1 \le p^{1+\epsilon} }} \ \sum_{\substack{C_2 \le c<2C_2\\|m_2|\le  \frac{ p^{2+\epsilon}}{N} \\  cm_2\equiv \alpha_{m_1}  \bmod m_1 }}  1  \ll p^\epsilon\Big(\frac{p}{q}+1\Big).
\]
This again is immediate, keeping in mind \eqref{crange}.

\bibliographystyle{amsplain}
\bibliography{twisted-subconvexity}

\providecommand{\bysame}{\leavevmode\hbox to3em{\hrulefill}\thinspace}
\providecommand{\MR}{\relax\ifhmode\unskip\space\fi MR }
\providecommand{\MRhref}[2]{%
  \href{http://www.ams.org/mathscinet-getitem?mr=#1}{#2}
}
\providecommand{\href}[2]{#2}
\begin{thebibliography}{10}

\bibitem{aggjonow}
K.~Aggarwal, Y.~Jo, and K.~Nowland, \emph{Hybrid level aspect subconvexity for
  {$GL(2)\times GL(1)$} {R}ankin-{S}elberg {$L$}-functions}, Hardy-Ramanujan J.
  \textbf{41} (2018), 104--117.

\bibitem{blohar}
V.~Blomer and G.~Harcos, \emph{Hybrid bounds for twisted {$L$}-functions}, J.
  Reine Angew. Math. \textbf{621} (2008), 53--79.

\bibitem{blohar2}
\bysame, \emph{Addendum: {H}ybrid bounds for twisted {$L$}-functions
  [mr2431250]}, J. Reine Angew. Math. \textbf{694} (2014), 241--244.

\bibitem{bloharmic}
V.~Blomer, G.~Harcos, and P.~Michel, \emph{A {B}urgess-like subconvex bound for
  twisted {$L$}-functions}, Forum Math. \textbf{19} (2007), no.~1, 61--105,
  Appendix 2 by Z. Mao.

\bibitem{bhkm}
V.~Blomer, P.~Humphries, R.~Khan, and M.~B. Milinovich, \emph{Motohashi's
  fourth moment identity for non-archimedean test functions and applications},
  Compos. Math. \textbf{156} (2020), no.~5, 1004--1038.

\bibitem{blokha2}
V.~Blomer and R.~Khan, \emph{Twisted moments of {$L$}-functions and spectral
  reciprocity}, Duke Math. J. \textbf{168} (2019), no.~6, 1109--1177.

\bibitem{blokha}
\bysame, \emph{Uniform subconvexity and symmetry breaking reciprocity}, J.
  Funct. Anal. \textbf{276} (2019), no.~7, 2315--2358.

\bibitem{bky}
V.~Blomer, R.~Khan, and M.~P. Young, \emph{Distribution of mass of holomorphic
  cusp forms}, Duke Math. J. \textbf{162} (2013), no.~14, 2609--2644.

\bibitem{blomil}
V.~Blomer and D.~Mili\'{c}evi\'{c}, \emph{{$p$}-adic analytic twists and strong
  subconvexity}, Ann. Sci. \'{E}c. Norm. Sup\'{e}r. (4) \textbf{48} (2015),
  no.~3, 561--605.

\bibitem{bur}
D.~A. Burgess, \emph{On character sums and {$L$}-series. {II}}, Proc. London
  Math. Soc. (3) \textbf{13} (1963), 524--536.

\bibitem{butkha}
J.~Buttcane and R.~Khan, \emph{{$L^4$}-norms of {H}ecke newforms of large
  level}, Math. Ann. \textbf{362} (2015), no.~3-4, 699--715.

\bibitem{byk}
V.~A. Bykovski\u{\i}, \emph{A trace formula for the scalar product of {H}ecke
  series and its applications}, Zap. Nauchn. Sem. S.-Peterburg. Otdel. Mat.
  Inst. Steklov. (POMI) \textbf{226} (1996).

\bibitem{coniwa}
J.~B. Conrey and H.~Iwaniec, \emph{The cubic moment of central values of
  automorphic {$L$}-functions}, Ann. of Math. (2) \textbf{151} (2000), no.~3,
  1175--1216.

\bibitem{dfi}
W.~Duke, J.~B. Friedlander, and H.~Iwaniec, \emph{Bounds for automorphic
  {$L$}-functions. {II}}, Invent. Math. \textbf{115} (1994), no.~2, 219--239.

\bibitem{hofloc}
J.~Hoffstein and P.~Lockhart, \emph{Coefficients of {M}aass forms and the
  {S}iegel zero}, Ann. of Math. (2) \textbf{140} (1994), no.~1, 161--181, With
  an appendix by Dorian Goldfeld, Hoffstein and Daniel Lieman.

\bibitem{houche}
F.~Hou and B.~Chen, \emph{Level aspect subconvexity for twisted
  {$L$}-functions}, J. Number Theory \textbf{203} (2019), 12--31.

\bibitem{ivi}
A.~Ivi\'{c}, \emph{On sums of {H}ecke series in short intervals}, J. Th\'{e}or.
  Nombres Bordeaux \textbf{13} (2001), no.~2, 453--468.

\bibitem{iwa}
H.~Iwaniec, \emph{Topics in classical automorphic forms}, Graduate Studies in
  Mathematics, vol.~17, American Mathematical Society, Providence, RI, 1997.

\bibitem{iwakow}
H.~Iwaniec and E.~Kowalski, \emph{Analytic number theory}, American
  Mathematical Society Colloquium Publications, vol.~53, American Mathematical
  Society, Providence, RI, 2004.

\bibitem{iwamic}
H.~Iwaniec and P.~Michel, \emph{The second moment of the symmetric square
  {$L$}-functions}, Ann. Acad. Sci. Fenn. Math. \textbf{26} (2001), no.~2,
  465--482.

\bibitem{kiryou}
E.~M. Kiral and M.~P. Young, \emph{The fifth moment of modular
  {$L$}-functions}, preprint (arXiv:1701.07507).

\bibitem{liumasyou}
S.-C. Liu, R.~Masri, and M.~P. Young, \emph{Rankin-{S}elberg {$L$}-functions
  and the reduction of {CM} elliptic curves}, Res. Math. Sci. \textbf{2}
  (2015), Art. 22, 23.

\bibitem{micven}
P.~Michel and A.~Venkatesh, \emph{The subconvexity problem for {${\rm GL}_2$}},
  Publ. Math. Inst. Hautes \'{E}tudes Sci. (2010), no.~111, 171--271.

\bibitem{mil}
D.~Mili\'{c}evi\'{c}, \emph{Sub-{W}eyl subconvexity for {D}irichlet
  {$L$}-functions to prime power moduli}, Compos. Math. \textbf{152} (2016),
  no.~4, 825--875.

\bibitem{munsin}
R.~Munshi and S.~K. Singh, \emph{Weyl bound for {$p$}-power twist of {$\rm
  GL(2)$} {$L$}-functions}, Algebra Number Theory \textbf{13} (2019), no.~6,
  1395--1413.

\bibitem{pen}
Z.~Peng, \emph{Zeros and central values of automorphic {L}-functions}, ProQuest
  LLC, Ann Arbor, MI, 2001, Thesis (Ph.D.)--Princeton University.

\bibitem{petyou2}
I.~Petrow and M.~P. Young, \emph{The fourth moment of dirichlet {$L$}-functions
  along a coset and the {W}eyl bound}, preprint (arXiv:1908.10346).

\bibitem{petyou1}
\bysame, \emph{The {W}eyl bound for {D}irichlet {$L$}-functions of cube-free
  conductor}, preprint (arXiv:1811.02452).

\bibitem{petyou3}
I.~Petrow and M.~P. Young, \emph{A generalized cubic moment and the {P}etersson
  formula for newforms}, Math. Ann. \textbf{373} (2019), no.~1-2, 287--353.

\bibitem{you}
M.~P. Young, \emph{Weyl-type hybrid subconvexity bounds for twisted
  {$L$}-functions and {H}eegner points on shrinking sets}, J. Eur. Math. Soc.
  (JEMS) \textbf{19} (2017), no.~5, 1545--1576.

\end{thebibliography}

\end{document}